\documentclass[a4paper,11pt]{article}
\usepackage{multirow}
\usepackage{color}
\usepackage{latexsym}
\usepackage{amssymb}
\usepackage{amsfonts}
\usepackage{amsmath}
\usepackage{amsthm}
\usepackage{indentfirst}
\usepackage{mathrsfs}
\usepackage{verbatim}

\newtheorem {theorem}{Theorem}[section]

\newtheorem {lemma}{Lemma}[section]
\newtheorem {proposition}{Proposition}[section]
\newtheorem {corollary}{Corollary}[section]

\theoremstyle{remark}

\usepackage[latin1]{inputenc}

\usepackage{fullpage}

\date{}

\author{Luca Ferrari\thanks{Dipartimento di Matematica e Informatica ``U. Dini", University of Firenze, Firenze, Italy\quad \texttt{luca.ferrari@unifi.it}. Partially supported by INdAM - GNCS project ``Strutture discrete in informatica: permutazioni, parking functions, linguaggi formali, ipergrafi".}\and Francesco Verciani\thanks{Institut f\"ur Mathematik, Universit\"at Kassel, 34132 Kassel, Germany, \, \texttt{francesco.verciani@uni-kassel.de}}}

\title{A new approach to Naples parking functions through complete parking preferences}

\begin{document}
	
\maketitle
	
\begin{abstract}
Naples parking functions were introduced as a generalization of classical parking functions, in which cars are allowed to park backwards, by checking up to a fixed number of previous spots, before proceeding forward as usual.

In this work we introduce the notion of a complete parking preference, through which we are able to give some information on the combinatorics of Naples parking functions. Roughly speaking, a \emph{complete parking preference} is a parking preference such that, for any index $j$, there are more cars with preference at least $j$ than spots available from $j$ onward. We provide a characterization of Naples parking functions in terms of certain complete subsequences of them. As a consequence of this result we derive a characterization of permutation-invariant Naples parking functions which turns out to be equivalent to the one given by (Carvalho et al., 2021), but using a totally different approach (and language).
\end{abstract}
	
\section{Introduction}

Parking functions are classical combinatorial objects, originally introduced by Konheim and Weiss \cite{KW} in a computer science context, but quickly gaining importance for themselves because of their appearance in numerous other contexts, such as hyperplane arrangements \cite{St}, representation theory \cite{KT} and combinatorial geometry \cite{AW}. More recently, there has been a renewed and growing interest in the combinatorics of parking functions, especially with regards to various generalizations that have been proposed. An exhaustive list of articles dealing with parking functions in the last decade is almost impossible to provide. We content ourselves to cite \cite{CCHJR}, a lovely paper where a large number of possible variations on the theme of parking functions is proposed, with an engaging and appealing style particularly suitable for undergraduates and probably even high school students. All articles cited by \cite{CCHJR}, as well as all articles citing \cite{CCHJR}, would deserve to be included in the references of our paper.

One of the most studied variations of parking functions are the so-called \emph{Naples parking functions}, in which cars are also allowed to check up to a certain number of previous spots before proceeding further. The enumerative combinatorics of Naples parking functions is not completely understood yet: in the literature there are recursive formulas for the number of $k$-Naples parking functions of length $n$ \cite{CHJLRRR}, but finding a closed form seems to be a rather difficult problem, as the only known explicit espression is rather involved \cite{CHJKRSV}. In our work we propose a new approach to the enumerative combinatorics of Naples parking functions, which is based on the notion of a complete parking preference. Given a parking preference $\alpha$ of length $n$, it may happen for a certain position $j$ that the number of available spots from $j$ onwards is strictly less than the number of cars whose preference is $j$ or more. This of course prevents $\alpha$ from being a parking function, but not a Naples parking function. Such positions are therefore critical in studying Naples parking functions. A complete parking preference is precisely one in which \emph{every} position (except the first, of course) is critical in the above sense. Using the notion of complete parking preference we are then able to provide characterizations of both generic and permutation-invariant Naples parking functions.

After introducing and studying properties of the excess function of a parking preference (Section \ref{e_f}) and of complete parking preferences (Section \ref{complete}), we describe our characterization results in Section \ref{characterization}, which are then summarized in Section \ref{further_work}, where we also state our projects for further work. 

\bigskip

In closing this introduction, we collect the main definitions and notations that we need in our work.

A \emph{parking preference} $\alpha =(a_1 ,a_2 ,\ldots ,a_n)$ of length $n$ is a $n$-tuple of positive integers to be interpreted as follows: for an ordered list $c_1 ,c_2 ,\ldots c_n$ of $n$ cars attempting to park on $n$ numbered spots in a street, $a_i$ denotes the preferred spot for the $i$-th car $c_i$. The set of all parking preferences of length $n$ will be denoted $PP_n$. The total number of cars having preference $i$ will be denoted with $|\alpha |_i$. The \emph{classical} parking rule requires that each car tries to park into its preferred spot and, if this is already occupied, then proceeds forward and parks into the first available spot, if any, otherwise it simply exits the street (thus failing to park). If all cars park successfully using the above rule, then we say that $\alpha$ is a \emph{parking function} of length $n$, and we denote the set of all those with $PF_n$.

As an example, the parking preference $\alpha =(3,1,3,5,2,4,2)\in PP_7$ is a parking function, such that, for instance, $|\alpha |_3 =2$ and $|\alpha |_6 =0$. During the parking process, we have that car $c_1$ parks in spot 3, car $c_2$ parks in spot 1, car $c_3$ parks in spot 4 (because it finds spot 3 already occupied), and so on.

An elegant characterization of parking functions says that $\alpha$ is a parking function if and only if $\sum_{i=1}^{j}|\alpha |_i \geq j$, for all $j\leq n$. The latter condition is equivalent to the fact that the nondecreassing rearrangement $(b_1 ,b_2 ,\ldots ,b_n)$ of $\alpha$ is such that $b_j \leq j$, for all $j\leq n$. As a consequence, we get that any rearrangement of a parking function is a parking function as well.

The classical parking rule can be modified into the \emph{$k$-Naples parking rule}, by allowing each car to check up to $k$ spots before its preferred one, and then proceed forward if it does not find an available place. A \emph{$k$-Naples parking function} is a parking preference such that all cars park successfully following the $k$-Naples parking rule, and the set of $k$-Naples parking functions of length $n$ will be denoted $PF_{n,k}$. We remark that, unlike classical parking functions, a rearrangement of a Naples parking function is not a Naples parking function in general.

\section{The excess function of a parking preference}\label{e_f}

Let $\alpha =(a_1 ,a_2 ,\ldots,a_n )\in PP_n$. Then, for every $j\leq n$, we set
\begin{equation}\label{excess}
u_\alpha (j)=\sum_{i=j}^{n}|\alpha|_i -(n-j+1).
\end{equation}

In other words, $u_\alpha (j)\in \mathbb{Z}$ denotes the excess of cars having preference at least $j$ with respect to the number of available spots from $j$ onwards. The map $u_{\alpha}$ will be called the \emph{excess function} of the parking preference $\alpha$. It is clear that, if $\alpha \in PF_n$, then $u_\alpha (j)\leq 0$ for all $j$ (and we will see in a moment that the reverse implication is also true). Thus those positions $j$ for which $u_\alpha (j)>0$ are somehow critical in trying to understand if the parking preference $\alpha$ is actually a Naples parking function. We set $U_{\alpha}=\{ j\leq n\, |\, u_\alpha (j)\geq 1\}$. In other words, $U_{\alpha}$ is the set of positions having positive excess function, and we remark that obviously $1\notin U_{\alpha}$, for all parking preferences $\alpha$.

\bigskip

In the next couple of lemmas we record some basic properties of the excess function that are frequently used in all the proofs of the results presented in this paper.

\begin{lemma}\label{elementary_properties}
	Let $\alpha \in PP_n$. For all $j,j_1 j_2 \leq n$, with $j_1 <j_2$, we have: 
	\begin{itemize}
		\item[(i)] $u_\alpha (j)=j-1-\sum_{i=1}^{j-1}|\alpha |_i$;
		\item[(ii)] $u_\alpha (j_2 )-u_\alpha (j_1 )=j_2 -j_1 -\sum_{i=j_1 }^{j_2 -1}|\alpha |_i$;
		\item[(iii)] $u_\alpha (j)=u_\alpha (j+1)+|\alpha |_j -1$, when $j<n$;
		\item[(iv)] $u_\alpha (j)<j$.
	\end{itemize} 
\end{lemma}

\begin{proof}
	Since of course $n=\sum_{i=1}^{j-1}|\alpha |_i +\sum_{i=j}^{n}|\alpha |_i$, using (\ref{excess}) we immediately get $(i)$. Moreover, again using the definition of $u_\alpha (j)$ for $j=j_2,j_1$, we get $(ii)$. Finally, we observe that $(iii)$ is a special case of $(ii)$ (when $j_1 =j$ and $j_2 =j+1$) and $(iv)$ clearly follows from $(i)$. 
\end{proof}

\begin{lemma}\label{elementary_intervals}
	Let $\alpha \in PP_n$ and let $[p,q]$ be a maximal interval of consecutive positions contained in $U_\alpha$. Then we have: 
	\begin{itemize}
		\item[(i)] $u_\alpha (p)=1$;
		\item[(ii)] $u_\alpha (p-1)=0$;
		\item[(iii)] $|\alpha |_{p-1}=0$;
		\item[(iv)] $|\alpha |_q \geq 2$.
	\end{itemize}     
\end{lemma}

\begin{proof}
	Recall that $1\notin U_{\alpha}$, hence $p\geq 2$. Since $[p,q]$ is maximal, we have that $p-1\notin U_{\alpha}$, and so $u_{\alpha}(p-1)\leq 0$ and $u_{\alpha}(p)\geq 1$. Using $(iii)$ of the previous lemma we thus get:
	\[
	1\leq u_{\alpha}(p)=u_{\alpha}(p-1)-|\alpha |_{p-1}+1\leq 1,
	\]
	which gives $(i)$. Moreover, we also have:
	\[
	0\geq u_{\alpha}(p-1)=u_{\alpha}(p)+|\alpha |_{p-1}-1=|\alpha |_{p-1}\geq 0,
	\]
	which gives $(ii)$ and $(iii)$.
	In order to prove $(iv)$, recall that $u_{\alpha}(q)\geq 1$. If $q=n$, then by definition of the excess function we have that $|\alpha |_n =u_{\alpha}(n)+1\geq 2$. If $q<n$, then clearly $u_{\alpha}(q+1)\leq 0$ and so $|\alpha |_q\geq |\alpha |_q +u_{\alpha}(q+1)=u_{\alpha}(q)+1\geq 2$, which gives $(iv)$.	
\end{proof}

Using $(i)$ of Lemma \ref{elementary_properties}, we get that the inequalities $\sum_{i=1}^{j}|\alpha |_i \geq j$ are equivalent to the inequalities $u_{\alpha}(j+1)\leq 0$, at least for all $j=1,\ldots n-1$. Moreover, since we know that $u_{\alpha}(1)=0$ and obviously $\sum_{i=1}^{n}|\alpha |_i =n$, we have that $\alpha$ is a parking function if and only $u_{\alpha}(j)\leq 0$ for all $j\leq n$. We record (and rephrase) this characterization of parking functions in the next proposition.

\begin{proposition}\label{easy_characterization}
	A parking preference $\alpha$ is a parking function if and only if $U_{\alpha}=\emptyset$.
\end{proposition}

We may now wonder whether a similar result holds also for Naples parking functions. Since Naples parking functions are \emph{not} in general permutation-invariant, a characterization in terms of the excess function only is not possible. However, we observe that the following necessary condition holds, which clearly generalizes Proposition \ref{easy_characterization}.

\begin{proposition}\label{necessary_PF}
	If $\alpha \in PF_{n,k}$, then $u_\alpha (j)\leq k$, for all $j\leq n$.
\end{proposition}

\begin{proof}
	By contradiction, suppose that there exists $j\leq n$ such that $u_\alpha (j)>k$. That is, the number of cars whose preference is at least $j$ is strictly greater than the number of available spots from $j$ onwards plus $k$. By $(i)$ of Lemma \ref{elementary_properties} we know that $u_\alpha (i)\leq i-1$, for all $i$, and so in particular $j\geq k+2$, i.e. there are at least $k+1$ parking spots available before the $j$-th one. Since $\alpha \in PF_{n,k}$, all cars must be able to park, and our hypothesis implies that there is at least one car (among those with preference $\geq j$) that needs to check at least up to $k+1$ spots before its preference. This is however impossible, since $\alpha$ is a $k$-Naples parking function.   
\end{proof}

Unfortunately, the above condition is not sufficient, as it can be realized by checking that the parking preference $(2,3,3)$ is not a 1-Naples parking function (and $u_\alpha (j)=1$, for $j=2,3$). Such a condition is however sufficient to show that there is at least one rearrangement of the parking preference that is a Naples parking function.

\begin{proposition}
	Suppose that $\alpha \in PP_n$ is such that $u_\alpha (j)\leq k$, for all $j\leq n$. If $\alpha$ is nonincreasing, then $\alpha \in PF_{n,k}$. 
\end{proposition}

\begin{proof}
	Suppose by contradiction that $\alpha \notin PF_{n,k}$, that is there is at least one index $i$ such that car $c_i$ is not able to park. This means that, when $c_i$ tries to park, all the spots in $[a_i -k,n]$ are already occupied. Since $\alpha$ is nonincreasing, all such spots have been occupied by cars with preference in $[a_i ,n]$. As a consequence, the number of cars with preference in $[a_i ,n]$ is at least $n-(a_i -k-1)+1=n-a_i +k+2$, hence
	\[
		u_{\alpha}(a_i )=\sum_{j=a_i}^{n}|\alpha|_{j} -(n-a_i +1)\geq (n-a_i +k+2)-(n-a_i +1)=k+1,
	\]
	which contradicts the hypothesis.
\end{proof}

\section{Complete parking preferences}\label{complete}

We are now ready to introduce complete parking preferences. These are a special subclass of parking preferences that play a crucial role in our characterization of Naples parking functions.

Given $n\geq 2$, a parking preference $\alpha \in PP_n$ is said to be \emph{complete} when $U_\alpha =[2,n]$. In other words, a complete parking preference is such that the set $U_\alpha$ is as large as possible (remember that $1\notin U_{\alpha}$). As an example, the parking preference $(5,3,3,5,4)\in PP_5$ is complete, whereas $(5,3,3,4,4)$ is not.

In order to collect some useful properties of complete parking functions, we need to introduce one more piece of notation. For a parking preference $\alpha =(a_1 ,a_2 ,\ldots ,a_n)\in PP_n$, and for a given integer $k\geq 0$, we consider the \emph{outcome map} $\psi^{(\alpha )}_k$, which is defined by setting $\psi^{(\alpha )}_k (i)$ to be the parking spot occupied by car $c_i$ following the $k$-Naples parking rule, with the convention that $\psi^{(\alpha )}_k (i)=\infty$ if car $c_i$ is unable to park. This is analogous to the usual outcome map defined for classical (i.e. $k=0$) parking preferences \cite{CDMY}, with a slight modification concerning the cars that are unable to park.

\bigskip

Our first goal is to provide a characterization of complete parking preferences that are also $k$-Naples parking functions. To this aim, the following lemma will be useful. It shows an interesting property of cars parking by driving forward, namely they create a sort of splitting of the set of spots.

\begin{lemma}\label{drive_forward}
	Let $\alpha \in PP_n$ and assume that all cars follow the $k$-Naples parking rule (for some $k\geq 1$). Let $i,j\leq n$ be such that $\psi^{(\alpha )}_k (i)=j$ and $a_i <j$ (i.e., car $c_i$ parks by driving forward). Then there are no cars with preference in $[j,n]$ that park in any spot in $[1,j]$. Moreover, if $\alpha \in PF_{n,k}$, then $u_{\alpha}(j)\leq -1$. 
\end{lemma}

\begin{proof}
	Let $h\leq n$ be such that $a_h \in [j,n]$. If $h<i$, then spot $j$ is still unoccupied, so car $c_h$ cannot park before position $j$. Now suppose that $h>i$. Since car $c_i$ has parked by driving forward, it has checked all the spots in the interval $[a_i -k,j]$, finding them already occupied. This implies that car $c_h$ cannot find any reachable free spots before position $j$ (since $a_h \geq j$ and $c_h$ can check up to $k$ spot before its preferred one).
	Finally, if $\alpha \in PF_{n,k}$, then all cars are able to park, and in particular all spots in $[1,j]$ are occupied by cars having preferences in $[1,j-1]$. Therefore there are at least $j$ cars having preferences in $[1,j-1]$ and so (using $(i)$ of Lemma \ref{elementary_properties})
	\[
		u_{\alpha}(j)=j-1-\sum_{i=1}^{j-1}|\alpha |_j \leq -1 ,
	\]
	as desired.    
\end{proof}

We now prove our characterization, which essentially states that complete Naples parking functions are those in which \emph{every} car parks by driving (weakly) backwards.

\begin{proposition}\label{char_complete}
	Let $\alpha \in PP_n$ be complete, and assume that all cars follow the $k$-Naples parking rule. The following are equivalent:
	\begin{itemize}
		\item[(i)] $\alpha \in PF_{n,k}$;
		\item[(ii)] for all $j\leq n$, spot $j$ is occupied by a car having preference at least $j$;
		\item[(iii)] for all $j\leq n$, $\psi^{(\alpha )}_k (j )\leq a_j$.
	\end{itemize}
\end{proposition}

\begin{proof}
	Suppose that $\alpha \in PF_{n,k}$. If there is some $j\leq n$ such that spot $j$ is occupied by a car whose preference is strictly less than $j$, then Lemma \ref{drive_forward} tells us that $u_{\alpha}(j)\leq -1$, contradicting the fact that $\alpha$ is complete. This proves that $(i)$ implies $(ii)$. The fact that $(ii)$ implies $(iii)$ is obvious (notice, in particular, that by $(ii)$ all spots are occupied, which excludes the possibility that $\psi^{(\alpha )}_k (j )=\infty$ for any $j$). Finally, if $(iii)$ holds, then necessarily $\psi^{(\alpha )}_k (j )<\infty$ for all $j$, which means that all cars are able to park, i.e. $\alpha \in PF_{n,k}$.     
\end{proof}

We are also able to refine the above characterization by providing more precise quantitative results about it.

\begin{proposition}\label{quantitative}
	Let $\alpha \in PP_n$ be complete, and assume that all cars follow the $k$-Naples parking rule. Then, for all $j\leq n$,
	\begin{equation}\label{back_complete}
	|\{ i\leq n\, |\, a_i \geq j, \psi^{(\alpha )}_k (i)<j \}|\leq u_\alpha (j).
	\end{equation}	
	
	In other words, the number of cars having preference at least $j$ which park in a spot strictly less than $j$ is at most $u_\alpha (j)$. 
	Moreover, if $\alpha \in PF_{n,k}$, then in (\ref{back_complete}) the equality holds.
\end{proposition}

\begin{proof}
	We start by observing that, if $j$ is a spot which remains unoccupied, then the left hand side of (\ref{back_complete}) is 0, and the thesis for $j$ immediately follows because $\alpha$ is complete. Therefore in the rest of the proof we assume that spot $j$ is occupied by some car.
	
	Suppose first that $j=n$. If spot $n$ is occupied by a car whose preference is strictly less than $n$, then Lemma \ref{drive_forward} implies that there are no cars with preference $n$ which are able to park, hence in particular the left hand side of (\ref{back_complete}) is 0 and the thesis follows for $n$. On the other hand, if spot $n$ is occupied by a car having preference $n$, then by definition of the excess function there are exactly $u_{\alpha}(n)$ more cars having preference $n$, i.e. $u_{\alpha}(n)$ is an upper bound for the number of cars having preference $n$ which are able to park somewhere strictly before $n$.
	
	Now let $j<n$ and suppose by induction that (\ref{back_complete}) holds for $j+1$. If spot $j$ is occupied by a car whose preference is strictly less than $j$, we can again invoke Lemma \ref{drive_forward} to get the thesis for $j$. Suppose instead that spot $j$ is occupied by a car having preference at least $j$. In this case we have:
	\begin{align*}
	|\{ i\leq n\, |\, a_i \geq j, \psi^{(\alpha )}_k (i)<j \}|&=|\{ i\leq n\, |\, a_i =j, \psi^{(\alpha )}_k (i)<j \}|+|\{ i\leq n\, |\, a_i \geq j+1, \psi^{(\alpha )}_k (i)<j \}| \\
	&\leq |\alpha |_j +|\{ i\leq n\, |\, a_i \geq j+1, \psi^{(\alpha )}_k (i)<j+1 \}| \\
	&\leq |\alpha |_j +u_{\alpha}(j+1),
	\end{align*}
	where the last inequality follows from the inductive hypothesis. Since we are assuming that spot $j$ is occupied by a car $c$ having preference at least $j$, there are two possible options: either $c$ has preference $j$, in which case $|\{ i\leq n\, |\, a_i =j, \psi^{(\alpha )}_k (i)<j \}|\leq |\alpha |_j -1$, or $c$ has preference strictly larger than $j$, in which case $|\{ i\leq n\, |\, a_i \geq j+1, \psi^{(\alpha )}_k (i)<j \}|<|\{ i\leq n\, |\, a_i \geq j+1, \psi^{(\alpha )}_k (i)<j+1 \}|$, and so $|\{ i\leq n\, |\, a_i \geq j+1, \psi^{(\alpha )}_k (i)<j \}|\leq  u_{\alpha}(j+1)-1$. We then have (also using $(iii)$ of Lemma \ref{elementary_properties})
	\[
		|\{ i\leq n\, |\, a_i \geq j, \psi^{(\alpha )}_k (i)<j \}|\leq |\alpha |_j +u_{\alpha}(j+1)-1=u_{\alpha}(j),
	\]
	as desired.
	
	To conclude, we observe that, if $\beta \in PF_{n,k}$, then all cars are able to park, and so, by definition of the excess function, $u_{\beta}(j)$ is a lower bound for the number of cars having preference at least $j$ and parking strictly before spot $j$. So, for a generic $k$-Naples parking function, the reverse of inequality (\ref{back_complete}) holds, hence for a complete $k$-Naples parking function the equality holds. 	 
\end{proof}

%
%

It is worth observing that, if $\alpha \in PP_n$ is not complete, then inequality (\ref{back_complete}) in general does not hold. Indeed, let $\alpha =(3,4,4,4,3)$ and suppose that $k=3$. Such a preference is not complete, since $u_{\alpha}(5)=-1$. Moreover, $u_{\alpha}(4)=1$, and there are two cars (namely $c_3$ and $c_4$) that park strictly before spot 4. Notice that, in this specific example, we also have that $\alpha$ is a $3$-Naples parking function, whose outcome map is $(3,4,2,1,5)$.

\section{A characterization of $k$-Naples parking functions through complete subsequences}\label{characterization}

The present section contains our main result, namely a characterization of Naples parking functions which makes use of the notion of complete parking preference that we have introduced and studied in the previous section. In order to prove such a characterization we need a preliminary result, which is of interest in its own. The following proposition gives a characterization of $k$-Naples parking functions in terms of the filling of certain critical spots which only depends on the sets $U_{\alpha}$.  

\begin{proposition}\label{p-1_spot}
		Given $\alpha \in PP_n$ and $k\geq 1$, we have that $\alpha$ is a $k$-Naples parking function if and only if, for every maximal interval $[p,q]\subseteq U_{\alpha}$, spot $p-1$ is occupied.
\end{proposition}

\begin{proof}
	Obviously, if $\alpha$ is a $k$-Naples parking function, then every spot is occupied by some car.
	
	On the other hand, suppose that $\alpha$ is not a $k$-Naples parking function. Let $h$ be the smallest spot which remains unoccupied. Then all spots in $[1,h-1]$ are occupied by cars having preferences in $[1,h-1]$, and there is no car having preference in $[1,h-1]$ which parks after spot $h$. We thus have (by $(i)$ of Lemma \ref{elementary_properties}) that $u_{\alpha}(h)=h-1-\sum_{i=1}^{h-1}|\alpha |_i=0$, i.e. $h\notin U_{\alpha}$. Moreover, since $h$ remains unoccupied, then $|\alpha |_h =0$, which gives (using $(iii)$ of Lemma \ref{elementary_properties}) $u_{\alpha}(h+1)=u_{\alpha}(h)-|\alpha |_h +1=1$, i.e. $h+1$ is the minimum of a (maximal) interval of $U_{\alpha}$.  
\end{proof}

Notice that, if $\alpha \in PF_{n,k}$, then all spots are occupied, and (by Lemma \ref{elementary_intervals}) $u_{\alpha}(p-1)=0$ and $|\alpha |_{p-1}=0$, hence Proposition \ref{drive_forward} implies that spot $p-1$ must be occupied by a car driving backwards (i.e., having preference at least $p$). 
The above proposition thus identifies a set of somehow ``critical spots", that need to be occupied by cars driving backwards in order for the parking preference to be a Naples parking function. This condition can be accomplished, for instance, by choosing a sufficiently large $k$ (or by ordering the cars in a suitable way, of course).

From the point of view of the previous proposition, classical parking functions are characterized by the fact that the above set of critical spots is empty, which is consistent with the fact that $U_{\alpha}=\emptyset$.   

\bigskip

The last thing we need before stating our characterization of Naples parking functions is a couple of technical definitions, followed by a technical (but very useful) lemma.

Given $\alpha =(a_1 ,a_2 ,\ldots ,a_n ) \in PP_n$ and $w\in \mathbb{Z}$, the \emph{$w$-shift} of $\alpha$ is the $n$-tuple $\tau_w (\alpha )=(a_1 -w, a_2 -w,\ldots ,a_n -w)$. This is particularly meaningful when $w$ is nonnegative and $w\leq \min (a_1 ,a_2 ,\ldots ,a_n )$. Moreover, it will be useful to extend the above notation to sets: given $W\subseteq \mathbb{Z}$, we set $\tau_w (W)=\{ i\in \mathbb{Z}\, |\, i+w\in W\}$. Finally, if $J=\{ j_1 ,j_2 ,\ldots j_h \} \subseteq [n]$, with $j_1 <j_2 <\cdots <j_h$, the \emph{restriction} of $\alpha$ to $J$ is the $h$-tuple $\alpha_{|_{J}}=(a_{j_1},a_{j_2},\ldots ,a_{j_h})$.

Now let $\alpha \in PP_n$ and suppose that $j\leq n$ is such that $u_{\alpha}(j)=0$. Set $J=\{ i\leq n\, |\, a_i\geq j\}$ and $J^c =[1,n]\setminus J$. Then it follows immediately from the definition of the excess function that $|J|=n-j+1$, and so $\alpha_{|_{J^c}}\in PP_{j-1}$ and $\tau_{j-1}(\alpha_{|_{J}})\in PP_{n-j+1}$. Moreover, for all $i$ for which it makes sense, an easy computation shows that $u_{\alpha_{|_{J^c}}}(i)=u_{\alpha}(i)$ (and so $U_{\alpha_{|_{J^c}}}=U_{\alpha}\cap [1,j-1]$) and $u_{\tau_{j-1}(\alpha_{|_{J}})}(i)=u_{\alpha}(i+j-1)$ (and so $U_{\tau_{j-1}(\alpha_{|_{J}})}=\tau_{j-1}(U_{\alpha})\cap [1,n-j+1]=
\tau_{j-1}(U_{\alpha}\cap [j,n])$). Therefore, roughly speaking, every parking preference $\alpha$ can be decomposed into two smaller parking preferences in correspondence of every position $j$ for which $u_{\alpha}(j)=0$, whose corresponding excess functions are essentially the same as the excess function of $\alpha$  (up to a suitable shift when needed). When $\alpha$ is a Naples parking function, one of these two preferences is guaranteed to be a Naples parking function as well.

\begin{lemma}\label{restricted_translated_pf}
	Let $\alpha \in PF_{n,k}$ and suppose that $j\leq n$ is such that $u_{\alpha}(j)=0$. Then (with the above notations) $\tau_{j-1}(\alpha_{|_{J}})\in PF_{n-j+1,k}$.
\end{lemma}

\begin{proof}
	When $j=1$ the lemma is obvious, so we can suppose that $j\geq 2$. Suppose by contradiction that $\tau_{j-1}(\alpha_{|_{J}})$ is \emph{not} a $k$-Naples parking function. Denote with $h$ the minimum spot which remains empty in the restricted problem (where of course $1\leq h\leq n-j+1$). We observe that, in the restricted problem, all cars having preferences in $[1,h-1]$ are able to park, and they necessarily park in the spots between 1 and $h-1$, which are all filled. Thus the number of cars having preferences between 1 and $h-1$ is precisely $h-1$, and consequently the number of cars having preference at least $h$ coincides with the number of available spots from $h$ to $n-j-1$. In other words, we have that $u_{\tau_{j-1}(\alpha_{|_{J}})}(h)=u_{\alpha}(h+j-1)=0$. Now, since $\alpha \in PF_{n,k}$, in the original problem the spot $h+j-1$ is occupied by some car. We observe that it cannot be occupied by a car having preference $<h+j-1$, otherwise, by Lemma \ref{drive_forward}, we would have $u_{\alpha}(h+j-1)\leq -1$, against what we have proved above. Therefore spot $h+j-1$ is occupied by some car driving backward, i.e. having preference $>h+j-1$ (notice that there are no cars with preference $h+j-1$, since we are supposing that in the restricted problem there are no cars having preference $h$). However, before spot $h$ is occupied in the original problem, the behavior of the cars having preference $>h+j-1$ is the same in the original and in the restricted problem. This implies that the car which occupies spot $h$ in the original problem will do the same in the restricted problem, against the hypothesis. 
\end{proof}

It is worth observing that the above lemma does \emph{not} imply that, if $u_{\alpha}(j)=0$, then all cars having preference in $[j,n]$ necessarily park in positions belonging to $[j,n]$. Consider, for instance, the parking preference $\alpha =(4,4,3,2,3)\in PF_{5,1}$. We have that $u_{\alpha}(4)=0$, and the set of cars having preference at least 4 is $\{ c_1 ,c_2\}$. Clearly $\tau_3 (\alpha_{|_{\{ 1,2\}}})=(1,1)$ is a 1-Naples parking function, in accordance with the lemma. However, during the parking process for $\alpha$, car $c_2$ parks in spot $3\notin \{ 4,5\}$.

We also remark that a similar result does not hold in general for the parking preference $\alpha_{|_{J^c}}$. The same preference $\alpha$ considered above serves as a counterexample. Indeed, recalling that $u_{\alpha}(4)=0$, we have that $\alpha_{|_{\{ 3,4,5\}}}=\{ 3,2,3\} \notin PF_{3,1}$.

\bigskip

The next theorem contains the announced characterization of Naples parking functions in terms of complete parking subsequences.

\begin{theorem}\label{main_characterization}
	Let $\alpha \in PP_n$. We have that $\alpha \in PF_{n,k}$ if and only if, for every maximal interval $[p,q]\subseteq U_{\alpha}$, there exists $J\subseteq [n]$ such that $a_j \in [p,p-2+|J|]$ (for all $j\in J$) and $\tau_{p-2}(\alpha_{|_{J}})\in PF_{|J|,k}$ is a complete $k$-Naples parking function.
\end{theorem}

\begin{proof}
	We start by supposing that $\alpha \in PP_n$ satisfies the condition in the statement of the theorem. Let $[p,q]\subseteq U_{\alpha}$ be a maximal interval, and suppose that there exists a subset $J=\{ j_1 ,j_2 ,\ldots ,j_{|J|}\}$ of $[n]$ as above. In the hope of improving the clarity of our argument, we introduce some dedicated notations for the restricted preference. First of all, we denote with $\tilde{\psi}_k$ the outcome map of $\tau_{p-2}(\alpha_{|_{J}})$. Moreover, we denote with $\tilde{c}_i=c_{j_i}$ the $i$-th car in the restricted problem, and with $\tilde{a}_i=a_{j_i}-(p-2)$ the corresponding preference. Since by hypothesis $\tau_{p-2}(\alpha_{|_{J}})\in PF_{|J|,k}$, we know that, in the restricted problem, all spots in $[1,|J|]$ are occupied. We now want to show (by induction) that, in the original problem, the corresponding spots in $[p-1,|J|+(p-2)]$ are also occupied. More formally, we want to show that, for all $i\in [1,|J|]$, all spots of the form $\tilde{\psi}_k (\tilde{c}_{i})+(p-2)$ are occupied in the original problem. 
	
	We start by considering $i=1$. Since car $\tilde{c}_{1}$ is the first car in the restricted problem, clearly $\tilde{\psi}_k (\tilde{c}_{1})=\tilde{a}_1 =a_{j_1}-(p-2)$. As a consequence, in the original problem either $\psi^{(\alpha )}_k (a_{j_1})=a_{j_1}$ or car $c_{j_1}$ finds spot $a_{j_1}$ already occupied. In either cases, after car $c_{j_1}$'s turn, spot $a_{j_1}=\tilde{\psi}_k (\tilde{c}_{1})+(p-2)$ is occupied in the original problem. Now suppose that, for all $i\in [1,h-1]$, after car $c_{j_i}$'s turn, spot $\tilde{\psi}_k (\tilde{c}_{i})+(p-2)$ is occupied in the original problem. If $\tilde{\psi}_k (\tilde{c}_{h})=\tilde{a}_{h}=a_{j_h}-(p-2)$, then we can argue as in the base case above. Otherwise, since $\tau_{p-2}(\alpha_{|_{J}})$ is complete, by Proposition \ref{char_complete} we have that $\tilde{\psi}_k (\tilde{c}_{h})=\tilde{a}_h -m$, for some $m\leq k$. This implies that, in the restricted problem, all spots in the interval $[\tilde{a}_h -(m-1), \tilde{a}_h]$ are of the form $\tilde{\psi}_k (\tilde{c}_{i})$, for some $i<h$. Thus, using the induction hypothesis, we get that in the original problem, when car $c_h$ tries to park, all spots in the interval $[a_{j_h}-(m-1),a_{j_h}]$ are already occupied. This means that car $c_h$ reaches spot $a_{j_h}-m=\tilde{a}_h -m+(p-2)=\tilde{\psi}_k (\tilde{c}_{h})+(p-2)$, which is then certainly occupied after $c_h$'s turn (either by $c_h$ itself or by some previous car).
	
	We have thus shown that, in the original problem, all spots in the interval $[p-1,|J|+(p-2)]$ are occupied, so in particular spot $p-1$ is occupied. Repeating the same argument for all maximal intervals $[p,q]\subseteq U_{\alpha}$, using Proposition \ref{p-1_spot}, we can conclude that $\alpha \in PF_{n,k}$.
	
	\medskip
	
	Now suppose that $\alpha \in PF_{n,k}$ and let $[p,q]\subseteq U_{\alpha}$ be a maximal interval. We start by considering the special case $p=2$. Set $M=\min\{ m\leq n\, |\, [1,m]\subseteq \psi^{(\alpha )}_k (\{ c_i \, |\, a_i \in [1,m] \})\}$. In other words, $M$ is the minimum of the set of all indices $m$ such that all spots in $[1,m]$ are occupied by cars having preferences in $[1,m]$. Observe that such a set is nonempty, since $\alpha \in PF_{n,k}$, and so $n$ belongs to the set, hence $M$ is well defined. Moreover, using $(iii)$ of Lemma \ref{elementary_intervals}, we get that $|\alpha|_1=0$, and so $M>1$. Let $J=\{j_1 ,j_2 ,\ldots ,j_M \}$ be the set of the indices of the cars parking in the spots of the interval $[1,M]$, with $j_1 <j_2 <\cdots <j_M$. Our goal is to show that the restriction $\alpha_{|_{J}}$ is a complete $k$-Naples parking function.
	
	First of all, it is obvious that $\alpha_{|_{J}}\in PF_{M,k}$. In order to prove that $\alpha_{|_{J}}$ is complete, we now proceed by ``reverse induction" to show that, for all $h\in [2,M]$, $u_{\alpha_{|_{J}}}(h)\geq 1$. For the case $h=M$, we observe that we cannot have $|\alpha_{|_{J}}|_M =0$, since otherwise we would have in particular that all spots in $[1,M-1]$ would be occupied by cars having preference in $[1,M-1]$, which contradicts the minimality of $M$. Similarly, we cannot have $|\alpha_{|_{J}}|_M =1$. In fact, assuming that $c_j$ is the only car having preference $M$, spot $M$ is occupied either by a car having preference $<M$, or by $c_j$; in the former case Lemma \ref{drive_forward} would imply that $c_j$ is not able to park (which is not the case, since $\alpha \in PF_{n,k}$), whereas in the latter case we would still get that all spots in $[1,M-1]$ would be occupied by cars having preference in $[1,M-1]$. We can thus conclude that $|\alpha_{|_{J}}|_M \geq 2$, hence $u_{\alpha_{|_{J}}}(M)\geq 1$.
	
	Now consider $h\in [2,M-1]$ and suppose that $u_{\alpha_{|_{J}}}(i)\geq 1$ for all $i\geq h+1$. Our goal is to prove that $u_{\alpha_{|_{J}}}(h)\geq 1$. By $(iii)$ of Lemma \ref{elementary_properties} we have that $u_{\alpha_{|_{J}}}(h)=u_{\alpha_{|_{J}}}(h+1)+|\alpha_{|_{J}}|_h -1$. If $|\alpha_{|_{J}}|_h \geq 1$, then clearly $u_{\alpha_{|_{J}}}(h)\geq u_{\alpha_{|_{J}}}(h+1)\geq 1$. So suppose that $|\alpha_{|_{J}}|_h =0$. Again by the minimality of $M$, there must be at least one car having preference in $[1,h-1]$ which parks in a spot $i\geq h$. If $i\geq h+1$, then by Lemma \ref{drive_forward} $u_{\alpha_{|_{J}}}(i)\leq -1$, which contradicts our induction hypothesis. If instead $i=h$, then again by Lemma \ref{drive_forward} we have that $u_{\alpha_{|_{J}}}(h)\leq -1$, and so $u_{\alpha_{|_{J}}}(h+1)=u_{\alpha_{|_{J}}}(h)-|\alpha_{|_{J}}|_h +1\leq 0$, which again is in contrast with the induction hypothesis.
	
	To conclude the proof we now turn to the case of a general $p$. Given a maximal interval $[p,q]\subseteq U_{\alpha}$, set $\hat{J}=\{ i\leq n\, |\, a_i \geq p-1 \}$. By Lemma \ref{restricted_translated_pf}, and recalling that $u_{\alpha}(p-1)=0$ (by $(ii)$ of Lemma \ref{elementary_properties}), we have that $\tau_{p-2}(\alpha_{|_{\hat{J}}})\in PF_{n-(p-2),k}$ and $[2,q-(p-2)]$ is a maximal interval of $U_{\tau_{p-2}(\alpha_{|_{\hat{J}}})}$. Therefore we can exploit what we have done in the case $p=2$ to get that there exists a nonempty set $J\subseteq \hat{J}$ such that $\tau_{p-2}(\alpha_{|_{J}})\in PF_{|J|,k}$ is a complete $k$-Naples parking function. Moreover, since the preferences of all cars in the restricted problem clearly belong to the interval $[2,|J|]$, we also have that $a_j \in [p,|J|+p-2]$ for all $j\in J$, which concludes the proof.   
\end{proof}

\begin{corollary}\label{J_ineq}
	Suppose that $\alpha \in PF_{n,k}$. Then, for any maximal interval $[p,q]$ of $U_{\alpha}$, a set $J\subseteq [n]$ satisfying the conditions of the above theorem is such that $|J|\geq q-p+2$.
\end{corollary}

\begin{proof}
	Suppose first that $p=2$. In the proof of the above theorem we found a set $J$ satisfying the required conditions such that $|J|=M$, where $M$ is the minimum of the set of all indices $m$ such that all spots in $[1,m]$ are occupied by cars having preference in $[1,m]$. We observe that $u_{\alpha}(M+1)\notin U_{\alpha}$, since, by $(i)$ of Lemma \ref{elementary_properties},
	\[
	u_{\alpha}(M+1)=M-\sum_{i=1}^{M}|\alpha |_i \leq 0.
	\]
	As a consequence, referring again to the notations of the above theorem, $|J|\geq q$, which is our thesis when $p=2$.
	
	The general case ($p>2$) follows from the case $p=2$ and the last part of the proof of the previous theorem, and is thus left to the reader.  
\end{proof}

\bigskip

To illustrate Theorem \ref{main_characterization}, consider the parking preference $\alpha =(8,4,7,1,6,8,7,5,10,1)\in PP_{10}$, for which $U_{\alpha}=[4,7]$. It is immediate to see that $\alpha \in PF_{10,2}$. Set $J=\{ 2,3,5,7,8\}$, and observe that, for all $j\in J$, $a_j \in [4,7]$. Now we have that $\tau_{2}(\alpha_{|_{J}})=(2,5,4,5,3)$, and it is easy to check that this is a complete 2-Naples parking function.    
 
We also observe that, for a given maximal interval of $U_\alpha$, the set $J$ of the above theorem is not necessarily unique. For instance, given $\alpha$ as above, it can be easily checked that other possible choices for $J$ are $J=\{ 1,2,3,6,7,8 \}$ or $J=\{ 1,2,5,6,7,8 \}$.

\bigskip

A remarkable corollary of the above theorem is a neat characterization of permutation-invariant Naples parking functions in terms of the cardinalities of the maximal intervals of $U_{\alpha}$. The following lemma will be useful in the proof of such a corollary.

\begin{lemma}
	Let $\alpha \in PF_{n,k}$ be a complete $k$-Naples parking function. Then $\psi^{(\alpha )}_k (n)=1$ and $a_n \leq k+1$.
\end{lemma}

\begin{proof}
	Let $\psi^{(\alpha )}_k (n)=h$. If $h>1$, then the set of spots $[1,h-1]$ is nonempty, and all such spots are occupied by cars whose preferences belongs to $[1,h-1]$ (since a car with preference $\geq h$ would check spot $h$ before checking previous spots). Therefore the number of cars having preference in $[1,h-1]$ is at least $h-1$, hence $u_{\alpha} (h)=h-1-\sum_{i=1}^{h-1}|\alpha |_i \leq 0$, against the fact that $\alpha$ is complete. Thus $h=1$, as desired. Moreover we get that $1=\psi^{(\alpha )}_k (n)\geq a_n -k$, i.e. $a_n \leq k+1$.
\end{proof}

\begin{corollary}\label{per_inv_k}
	Let $\alpha =(a_1 ,a_2 ,\ldots ,a_n ) \in PP_n$ and $k\geq 1$. The following are equivalent:
	\begin{itemize}
		\item[(i)] for every maximal interval $[p,q]\subseteq U_{\alpha}$, $|[p,q]|\leq k$;
		\item[(ii)] for all permutations $\sigma$ of size $n$, the $n$-tuple $\sigma (\alpha )=(a_{\sigma (1)} ,a_{\sigma (2)} ,\ldots ,a_{\sigma (n)})$ is a k-Naples parking function of length $n$. In particular, $\alpha \in PF_{n,k}$.
	\end{itemize}
\end{corollary}

\begin{proof}
	We start by supposing that condition $(i)$ holds, and let $[p,q]$ be a maximal interval of $U_{\alpha}$. Then $u_{\alpha}(q+1)\leq 0$, and so by $(iii)$ of Lemma \ref{elementary_properties} we have that $|\alpha |_q =u_{\alpha}(q)-u_{\alpha}(q+1)+1\geq u_{\alpha}(q)+1$. Now choose arbitrarily $u_{\alpha}(q)+1$ cars among those having preference $q$ and let $H$ be the set of their indices. Set $J=\{ i\leq n\, |\, a_i \in [p,q-1]\} \cup H$. We want to show that $J$ satisfies the properties stated in Theorem \ref{main_characterization}. First of all we observe that, by $(i)$ of Lemma \ref{elementary_intervals}, $u_{\alpha}(p)=1$, and so, by $(ii)$ of Lemma \ref{elementary_properties}, 
	\[
	|J|=\sum_{i=p}^{q-1}|\alpha |_i +u_{\alpha}(q)+1=(q-p)+u_{\alpha}(p)+1=q-p+2.
	\]
	
	Now consider the parking preference $\tau_{p-2}(\alpha_{|_J})$. It is easy to see that, for all $j\in [1,q-p+2]$, $u_{\tau_{p-2}(\alpha_{|_J})}(j)=u_{\alpha}(j+(p-2))$. Since $u_{\alpha}(i)\geq 1$ for all $i\in [p,q]$, we thus have that $u_{\tau_{p-2}(\alpha_{|_J})}(j)\geq 1$ for all $j\in [2,q-p+2]$, i.e. $\tau_{p-2}(\alpha_{|_J})$ is complete. Moreover, since $|J|=q-p+2=|[p,q]|+1\leq k+1$, it is obvious that $\tau_{p-2}(\alpha_{|_J})$ is a $k$-Naples parking function. We can thus invoke Theorem \ref{main_characterization} to conclude that $\alpha \in PF_{n,k}$. Finally, we notice that the definition of $J$ does not depend on the specific order of the cars. In other words, for any permutation $\sigma \in S_n$, we can construct exactly the same set $J$ also for $\sigma (\alpha )=(a_{\sigma (1)} ,a_{\sigma (2)} ,\ldots ,a_{\sigma (n)})$, so also $\sigma (\alpha )$ is a $k$-Naples parking function.
	
	\medskip
	
	Conversely, suppose that there exists a maximal interval $[p,q]$ of $U_{\alpha}$ such that $|[p,q]|=q-p+1>k$. Since to disprove condition $(ii)$ we just need to show that there is at least one rearrangement of $\alpha$ which is not a $k$-Naples parking function, we can assume without loss of generality that $\alpha$ is weakly increasing. Suppose by contradiction that $\alpha \in PF_{n,k}$. By Theorem \ref{main_characterization} we thus have that there exists some $J$ such that $\tau_{p-2}(\alpha_{|_J})$ is a complete $k$-Naples parking function. Clearly also $\tau_{p-2}(\alpha_{|_J})$ is weakly increasing, and it is immediate to observe that the preference of the last car is equal to $|J|$ (indeed there must be at least one car having preference $|J|$, because $\tau_{p-2}(\alpha_{|_J})$ is complete). However, by the previous lemma, we should have $|J|\leq k+1$, whereas, by Corollary \ref{J_ineq}, we have $|J|\geq q-p+2>k+1$, hence the contradiction.      
\end{proof}

An equivalent characterization of permutation-invariant Naples parking functions  is given in \cite{CHRTV}, where the authors prove that a $k$-Naples parking function is permutation-invariant if and only if its ascending rearrangement is a $k$-Naples parking function. It is then shown that ascending parking preferences are related to Grand-Dyck paths, and in particular an ascending parking preference $\alpha$ is a $k$-Naples parking function if and only if in the corresponding Grand-Dyck path the length of every factor lying entirely below the $x$-axis (except at most for the starting and ending points) is at most $2k$. This is in fact equivalent to our characterization in terms of the cardinalities of the maximal intervals of $U_{\alpha}$, which has the advantage of being intrinsic to the parking preference $\alpha$ and does not involve other combinatorial structures. Moreover, the characterization given in Corollary \ref{per_inv_k} does not require to explicitly perform any parking process, but it just requires to check some counting property directly on the parking preference.

\section{Conclusions and further work}\label{further_work}

In the present work we have introduced the notion of \emph{complete parking preference}, which then has been used to provide a characterization of Naples parking functions, stated in Theorem \ref{main_characterization}. Together with Proposition \ref{p-1_spot}, they propose a method to check if a given parking preference is indeed a Naples parking function which is totally intrinsic to the preference, and only requires checking if it has some specific properties. We can summarize the results obtained in Theorem \ref{main_characterization}, Proposition \ref{p-1_spot} and Corollary \ref{per_inv_k} as follows.

\begin{theorem}
	Let $\alpha \in PP_n$ and suppose that all cars follow the $k$-Naples parking rule, with $k\geq 1$. Let $[p,q]\subseteq U_{\alpha}$ be a maximal interval. The following are equivalent:
	\begin{itemize}
		\item[(i)] restricting to cars having preferences at least $p$, spot $p-1$ is occupied by some car (by driving backwards, of course);
		\item[(ii)] there esists $J\subseteq [p,q]$, with $|J|\geq q-p+2$, such that $\tau_{p-2}(\alpha_{|_{J}})\in PF_{|J|,k}$ is a complete $k$-Naples parking function.
	\end{itemize} 
	Each of the above stated conditions is guaranteed to be satisfied when $|[p,q]|\leq k$. Moreover, $\alpha \in PF_{n,k}$ if and only if all maximal intervals of $U_{\alpha}$ having size at least $k+1$ satisfy any of the above conditions.
\end{theorem}

From the enumerative point of view, once complete parking preferences have been introduced, a natural task is to count complete $k$-Naples parking functions (with respect to the length). Moreover, since the notion of excess function of a parking preference is useful to characterize permutation-invariant Naples parking functions, we expect to be able to give enumerative results also concerning this last class of Naples parking function, thus enhancing and improving the results already achieved in \cite{CHRTV}. We plan to address these questions in a forthcoming paper.

\bigskip

To conclude, we would like to propose a far reaching generalization of Naples parking functions that in our opinion deserves to be investigated.

The classical $k$-Naples parking rule requires that \emph{all} cars check up to $k$ spots preceding the preferred one before proceeding forward. However, one may ask that \emph{each} car $c_i$ has to check up to a certain number $k_i$ (depending on $i$) of spots by driving backwards. In other words, we can consider a parking rule in which each car has its own number of spots to check by driving backwards. The resulting notion of parking function seems not to have been previously considered in the literature, and we also plan to define and study it in a forthcoming paper.

\end{document}